\documentclass[journal]{IEEEtran}

\usepackage[tbtags]{amsmath}
\usepackage{amsfonts}
\usepackage{eurosym}
\usepackage{graphics}
\usepackage{graphicx}
\usepackage{array}
\usepackage{shortvrb}
\usepackage{epsf}
\usepackage{graphicx}
\usepackage{rotating}
\usepackage{float}
\usepackage{color}
\usepackage{multirow}
\usepackage{float}
\usepackage{hhline}
\newcommand{\bfg}[1]{\mbox{\boldmath $#1$\unboldmath}}
\newcommand{\fraca}[2]{\displaystyle\frac{#1}{#2}}
\newcommand{\mm}[3]{\renewcommand{\arraystretch}{0.8}\begin{array}[t]{c}\mbox{#1}
\\ #2\end{array}\begin{array}[t]{c}#3\end{array}
\renewcommand{\arraystretch}{1}}
\def \R {{\rm I\kern -2.2pt R\hskip 1pt}}

\begin{document}
\title{Robust Dynamic Transmission and Renewable Generation Expansion Planning: {W}alking Towards Sustainable Systems}

\author{Cristina Rold\'an, Agust\'in A. S\'anchez de la Nieta,~\IEEEmembership{Member,~IEEE}, Roberto M\'inguez,
        and Raquel Garc\'ia-Bertrand,~\IEEEmembership{Senior Member,~IEEE}
\thanks{
\newline \indent
C. Rold\'an, A. A. S\'anchez de la Nieta, and R. Garc\'{\i}a-Bertrand are with the Department of Electrical Engineering, Universidad de Castilla-La
Mancha, Ciudad Real, Spain (e-mail: \mbox{Cristina.Roldan@uclm.es}, \mbox{agustinsnl@gmail.com}, \mbox{Raquel.Garcia@uclm.es})
\newline \indent
R. M\'{\i}nguez is with Hidralab Ingeniería y Desarrollo, S.L., Spin-Off UCLM, Hydraulics Laboratory, Universidad
 de Castilla-La Mancha, Ciudad Real E-13071, Spain (e-mail: \mbox{roberto.minguez@hidralab.com}).}
}

\markboth{}%
{}

\maketitle

\begin{abstract}
Nowadays, the transition from a conventional generation system to a renewable generation system is one of the most difficult challenges for system operators and companies. There are several reasons: the long-standing impact of investment decisions, the proper integration of renewable sources into the system, the present and future uncertainties, and the convenience to consider an integrated year-by-year representation of both uncertainties and investment decisions.
However, recent breakthroughs in Dynamic Transmission Network Expansion Planning (DTNEP) have demonstrated that the use of robust optimization might render this problem computationally tractable for real systems. This paper intends to consider not only the capacity expansion of lines, but the construction and/or dismantling of renewable and conventional generation facilities as well.
 The Dynamic Transmission Network and Renewable Generation Expansion Planning (DTNRGEP) problem is formulated as an adaptive robust optimization problem with three levels. First level minimizes the investment costs of transmission network and generation expansion planning, the second level maximizes system operational costs with respect to uncertain parameters, while the third level minimizes those operational costs with respect to operational decisions.
  The method is tested for two cases: i) an illustrative example based on Garver IEEE system and ii) a case study using the IEEE 118-bus system. Numerical results from these examples demonstrate that the proposed model allows making optimal decisions towards reaching a sustainable power system, while overcoming problem size limitations and computational intractability for realistic cases.
\end{abstract}


\vspace{-0.2cm}
\section*{Nomenclature}
\renewcommand{\labelitemi}{}
This section states the main notation used in this paper for quick
reference. 

{\bf Indices and Sets:}
\begin{description}
%
%
%
\item[$\mathcal{D}$] Set of indices of demand.

\item[$\mathcal{G}$] Set of indices of all generation units installed at the beginning of time horizon considered which can not be removed from the system.

\item[$\mathcal{G}^{+}$] Set of all prospective and independent new possible generators.

\item[$\mathcal{G}_{g}^{+}$] Set of all prospective new generators which can be installed at different phases associated with group $g$.

\item[$g$] Index for groups of generators built per phases.

\item[$\mathcal{G}^{-}$] Set of all generators to be uninstalled or dismantled during the study period.

\item [$i$] Index related to generators.

\item [$j$] Index associated with loads.

\item [$k$] Index referring to lines.

\item [$l$] Counter index for each iteration.

\item[$\mathcal{L}$] Set of all existing transmission lines at the beginning of time horizon considered.

\item[$\mathcal{L^{+}}$] Set of all prospective transmission lines.

\item[$\mathcal{N}$] Set of all networks buses.

\item[$n$] Index related to buses.

\item[$n(i)$] Bus index where the $i$-th generating unit is located.

\item[$n(j)$] Bus index where the $j$-th demand is located.

\item[$\mathcal{T}$] Set of indices of years.

\item[$\Psi_{n}^{\rm D}$] Set of indices of the demand located at bus $n$.

\item[$\Psi_{n}^{\rm G}$] Set of indices of the generating units located at bus $n$.

\item[$\mathcal{U}^{(t)}$] Set of indices of the uncertain variables for time period $t$.

\end{description}
{\bf Constants:}

\begin{description}

\item[$b_{k}$] Susceptance of line $k$ (S).

\item[$c^{\rm G}_{i}$] Generator $i$ operational cost (\euro/MWh).
\vspace{0.02cm}
\item[$c^{\rm GI}_{i}$] Generator $i$ investment cost (\euro).

\item[$c^{\rm LI}_{k}$] Line $k$ investment cost (\euro).
\vspace{0.02cm}
\item[$c^{\rm S}_{j}$] Consumer $j$ load-shedding cost (\euro/MWh).

\item[$e^{(t)}_j$] Percentage of load shed by the $j$-th demand for year $t$.

\item[$h_{\mu,j}^{(t)}$] Nominal value evolution factor  for demand $j$ and period $t$.

\item[$h_{\sigma,j}^{(t)}$] Dispersion value evolution factor for demand $j$ and period $t$.

\item[$f_{k}^{\rm max}$] Line $k$ capacity (MW).

\item[$I$] Discount rate.

\item[$N_{g}$] Total number of generators built per phases in the group $g$.

\item[$N_y$] Number of study periods.

\item[$o(k)$] Line $k$ sending-end bus.

\item[$r(k)$] Line $k$ receiving-end bus.

\item[$\Pi_{\rm L}$] Transmission expansion investment budget (\euro).

\item[$\Pi_{\rm G}$] Generation expansion investment budget (\euro).

\item[$\sigma$]  Annual weighting factor (h).

\item[$t_i^{G^{-}}$] Time period when generator $i\in \mathcal{G^{-}}$ is uninstalled or dismantled.

\end{description}

{\bf Primal variables:}
\begin{description}
\item[${\bfg u}^{(t)}$] Vector of random or uncertain parameters for year $t$, including generation capacities and loads (MW).

\item[$f^{(t)}_{k}$] Line $k$ power flow for year $t$ (MW).
\vspace{0.02cm}
\item[$g^{(t)}_{i}$] Power production of generating unit $i$ for year $t$ (MW).

\item[$p^{(t)}_{j}$] Power consumption of demand $j$ for year $t$ (MW).

\item[$r^{(t)}_{j}$] Load shed of demand $j$ for year $t$ (MW).

\item[$x_k^{(t)}$] Binary variable representing new line $k$ construction at the beginning of year $t$.

\item[$\tilde{x}_k^{(t)}$] Line $k$ status ({\em existing} vs {\em no existing}) at the beginning of year $t$.

\item[$y_k^{(t)}$] Binary variable representing new generator $i$ construction at the beginning of year $t$.

\item[$\tilde{y}_k^{(t)}$] Generator $i$ status ({\em existing} vs {\em no existing}) at the beginning of year $t$.

\item[$\theta^{(t)}_{n}$] Bus $n$ voltage angle for year $t$ (radians).

\end{description}

\section{Introduction}

\subsection{Motivation}
The new objective of Kyoto Protocol for reducing Greenhouse Gases (GHG) encourages the development of renewable energy sources within electric systems \cite{kyoto:1}. The main reason is to fight against the growing trend of worldwide average temperature and climate change, and thus, it is expected that vast amounts of new generation facilities, specially renewable, will be built in the medium-term future.

Transmission network and renewable generation expansion planning analyze the issue of how to expand or reinforce an existing power transmission network, incorporate new renewable generation facilities and dismantling the old ones to adequately service system loads over a given time horizon while decreasing GHG emissions. This problem is challenging for several reasons \cite{LumbrerasR:16}:
\begin{enumerate}
    \item Transmission and generation investment decisions have a long-standing impact on the power system as a whole.

    \item Transmission and generation investments, specially new generation sources, must be integrated appropriately into the existing system.

    \item Consumption and renewable energy generation uncertainties, such as with wind and solar power plants, complicate the problem resolution. Note that wind power has been the most developed renewable technology in the last decade, while the next renewable technology, in constant evolution, is the photovoltaic power. The introduction of this type of renewable sources in the generation mix increases the uncertainty about the feasibility of generation. Some references attempt to reduce the uncertainty by means of synergies between renewable generators, i.e. combined wind and hydro-pump generators \cite{de2013optimal}, or using physical bilateral contracts \cite{de2016optimal}.

    \item The expansion planning problem is by nature a multi-stage problem that entails planning a horizon of several years.  Keeping the full dynamic complexity of the problem has been considered to be highly complex, mostly resulting in computationally intractable problems.
\end{enumerate}
%

\subsection{Literature Review}
Transmission and generation expansion planning have been extensively studied areas from the beginning of power systems~\cite{wang1994modern}. These problems have been carried out by means of several mathematical programming techniques, such as stochastic programming \cite{birge2011introduction}, genetic algorithms \cite{whitley1994genetic}, or Adaptive Robust Optimization (ARO) \cite{ThieleTE:2}, among others.

Transmission expansion planning was first described in 1970 by \cite{Garver:70}. A new strategy for transmission expansion in a competitive electricity market is proposed in \cite{fang2003new}, while \cite{WuCX:08} studies the minimum load cutting problem existing in the process of transmission network expansion planning when load is uncertain within given intervals. A bi-level approach for transmission expansion planning within a market environment is proposed in \cite{GarcesCGR:09}. In an attempt to improve computational tractability, several robust approaches are presented in \cite{YuCW:11,Jabr:13,RuizC:15,MinguezG:15} by using ARO. They proved that computational tractability is possible for real-size systems. However, all these approaches take a wide spectrum of simplifying assumptions by considering static and sequential static models in an attempt to consider the year-by-year representation of investment decisions. In contrast, an alternative adaptive robust transmission network expansion planning formulation is proposed by \cite{MinguezGandBertrand:16}, which keeps the full dynamic complexity of the problem and reaches the global optimal solution of the problem.

Regarding generation expansion planning, it was first evaluated in 1981 by \cite{anders1981genetration}. Reference \cite{chen2004generation} attempted to solve the problem by using Lagrangian relaxation and probabilistic production simulation. An heuristic approach for power generation expansion planning with emission control is introduced in \cite{sirikum2006power}, while a modified version based on chance-constrained optimization is implemented by \cite{mazadi2009modified}. In contrast,  in \cite{nanduri2009generation} a two-level game-theoretic model is used. There are also robust approaches dealing with capacity expansion planning, a two-stage robust optimization model considering the uncertainty in investment costs is developed in \cite{dehghan2014two}.

Finally, joint consideration of transmission and generation expansion planning under uncertainty is checked in \cite{gorenstin1993power} and a comprehensive review is presented in \cite{hemmati2013comprehensive}. Transmission and generation expansion planning under risk using stochastic programming is provided by \cite{lopez2007generation}, and a market-based model with uncertainties is presented by \cite{roh2009market}. A static (one period) robust approach is modelled in \cite{dehghan2016reliability}, where electric demand and wind power generation are continuous uncertain variables, modelled through bounded intervals. In this case, availability of units and lines is represented by means of uncertain discrete variables, whose behaviour is modeled thorough probability distributions.

\subsection{Aims and Contributions}
The purpose of this paper is three-fold:
\begin{enumerate}
 \item To propose an effective model, which allows making optimal decisions in an integrated way to transform any power system into a sustainable system. It would allow to make decisions about where and when new transmission lines and/or new generation facilities, such as,  wind and photovoltaic power plants, should be installed during a long-term time frame of several years.
  \item To extend the ARO formulation associated with the dynamic expansion planning problem proposed by \cite{MinguezGandBertrand:16} for considering also the construction and/or dismantling of renewable and conventional generation facilities.
  \item To show that computational tractability for a year-by-year, multi-year or multi-stage representation of investment decisions ({\it dynamic approach}) associated with transmission and generation capacities is possible for real-size systems. In addition, global optimality is guaranteed.
\end{enumerate}

It is worth stressing that the proposed model is highly flexible in what regards to generation capacity expansion possibilities. It is possible to dismantle conventional generation facilities reaching their lifetimes during the time period considered, it is possible to consider the construction of renewable generation facilities in different phases, and the inclusion of conventional facilities without uncertainties. In summary, the model allows to take into consideration all practical aspects required to effectively transform any power system into a sustainable power system.

\subsection{Paper Structure}
The rest of the paper is structured as follows. Section~\ref{Model} describes the robust formulation of the DTNRGEP problem. The proposed decomposition method to solve the problem is described in Section~\ref{Descomposicion}. Section~\ref{s4} provides numerical results for two examples. Finally, in Section~\ref{s5} relevant conclusions are drawn.

\section{Robust Dynamic Transmission Network and Renewable Generation Expansion Planning Formulation}\label{Model}
A detailed formulation of the robust DTNRGEP problem is as follows:
\begin{align}
&\mm{Minimize}{x_k^{(t)},\tilde{x}_k^{(t)}, y_i^{(t)}\tilde{y}_i^{(t)}}\displaystyle\sum_{t \in \mathcal{T}} \fraca{1}{(1+I)^{t-1}}\left(\displaystyle\sum_{k \in \mathcal{L}^+}c_k^{\rm LI} x_k^{(t)}\right. \nonumber \\
&\left.\hspace{1cm}+\sum_{i \in \mathcal{G}^+\cup \mathcal{G}_{g}^{+};\forall g}c_i^{\rm GI} y_i^{(t)}+\fraca{c_{\rm op}^{(t)}}{(1+I)}\right); \label{eq1}
\end{align}
subject to
\begin{align}
\Pi_{\rm L} & \geq  \displaystyle\sum_{t \in \mathcal{T}} \displaystyle\sum_{k \in \mathcal{L}^+} \fraca{1}{(1+I)^{t-1}}c^{\rm LI}_k x_k^{(t)}  \label{eq2}\\
\tilde{x}_k^{(t)} &= 1;\;\forall k \in \mathcal{L},\forall t \in \mathcal{T} \label{eq2a}\\
\tilde{x}_k^{(t)} &= \sum_{p=1}^{p=t} x_k^{(p)};\;\forall k \in \mathcal{L}^+,\forall t \in \mathcal{T}\label{eq2b}\\
\sum_{t \in \mathcal{T}} x_k^{(t)} &\leq 1;\;\forall k \in \mathcal{L}^+\label{eq2c}\\
x_k^{(t)}  & \in  \{0,1\};\;\forall k \in \mathcal{L}^+,\forall t \in \mathcal{T} \label{eq2d}\\
\Pi_{\rm G} & \geq  \displaystyle\sum_{t \in \mathcal{T}} \displaystyle\sum_{i \in \mathcal{G}^+\cup \mathcal{G}_{g}^{+};\forall g} \fraca{1}{(1+I)^{t-1}}c^{\rm GI}_i y_i^{(t)}  \label{eq2f}\\
\tilde{y}_i^{(t)} &= 1;\;\forall i \in \mathcal{G},\forall t \in \mathcal{T} \label{eq2g}
\end{align}
\begin{align}
\tilde{y}_i^{(t)} &= \sum_{p=1}^{p=t} y_i^{(p)};\;\forall i \in \mathcal{G}^+\cup \mathcal{G}_{g}^{+};\forall g;\forall t \in \mathcal{T};\label{eq2h}\\
\sum_{t \in \mathcal{T}} y_i^{(t)} &\leq 1;\;\forall i \in \mathcal{G}^+\cup \mathcal{G}_{g}^{+};\forall g\label{eq2i}\\
y_i^{(t)}  & \in  \{0,1\};\;\forall i \in \mathcal{G}^+\cup \mathcal{G}_{g}^{+};\forall g;\forall t \in \mathcal{T} \label{eq2j}\\
\tilde{y}_i^{(t)} &= 1;\;\forall i \in \mathcal{G^-},\forall t=1,...,t_i^{G^{-}} \label{eq2k}\\
\tilde{y}_i^{(t)} &= 0;\;\forall i \in \mathcal{G^-},\forall t=t_i^{G^{-}}+1, ..., N_y \label{eq2l}\\
y_{i+1}^{(t)} &\leq \tilde{y}_{i}^{(t)};\forall i \in \mathcal{G}_{g}^{+};\forall g;\forall t \in \mathcal{T}  \label{eq2m}\\
y_{i+1}^{(t)}+y_{i}^{(t)} &\leq 1;\forall i \in \mathcal{G}_{g}^{+};\forall g;\forall t \in \mathcal{T}  \label{eq2n}
\end{align}
where the objective function (\ref{eq1}) is the sum of the present values of cost over the time horizon, i.e. net present cost (NPC).
Equations (\ref{eq2})-(\ref{eq2d}) are the constraints related to the construction of lines as presented in~\cite{MinguezGandBertrand:16}, which: i) limit the maximum expansion investment (\ref{eq2}), ii) force the line status to 1 for all existing transmission lines at the beginning (\ref{eq2a}), and iii) once the line has been constructed (\ref{eq2b}), iv) ensure that no line is constructed more than once throughout the time horizon considered (\ref{eq2c}), and v) establish the binary nature of line investment decisions (\ref{eq2d}).
Constraints (\ref{eq2f})-(\ref{eq2n}) are novel and associated with generation facilities. Constraint (\ref{eq2f}) keeps the maximum amount of generation investment within the available budget. Constraints (\ref{eq2g}) and (\ref{eq2h}) make the generation status  equal to 1 for all existing generation facilities at the beginning of the time horizon considered which will not be dismantled, and once the generation facility has been constructed, respectively, while constraint (\ref{eq2i}) ensures that no generation facility is constructed more than once. Constraint (\ref{eq2j}) establishes the binary nature of generation investment decisions. For generators to be dismantled during the study period ($\forall i \in \mathcal{G^-}$), constraint (\ref{eq2k}) make the generation $i$ status equal to 1 until the facility is dismantled, i.e. $t\le t_i^{G^{-}}$, while constraint (\ref{eq2l}) makes the status equal to 0 once it is dismantled, i.e. $t> t_i^{G^{-}}$. Constraints (\ref{eq2m}) and (\ref{eq2n}) ensure that for each generation group to be constructed in consecutive phases, the order of construction is sequential according to the generator set $\mathcal{G}_{g}^{+}$ order, thus phase $i+1$ can not be constructed before phase $i$.
Note that it would be straightforward to consider the possibility to dismantle old lines analogously to the generation case by adapting constraints (\ref{eq2k}) and (\ref{eq2l}).

Given the values of the first-stage decision variables $\tilde{x}_k^{(t)}, \tilde{y}_i^{(t)}$, operational costs $c_{\rm op}^{(t)}$ in (\ref{eq1}) for each period $t;\;\forall t \in \mathcal{T}$ are obtained using the following optimization problem. Note that the dual variables associated with constraints are provided separated by a colon.
\begin{align}
c_{\rm op}^{(t)}=\!\!\!\!\mm{Maximum}{{u}^{(t)}\in \mathcal{U}^{(t)}}\!\!\!\!\!\!\!\!&\mm{Minimum}{g^{(t)}_i, p^{(t)}_j, r^{(t)}_j, \\ \theta^{(t)}_n, f^{(t)}_k}\!\!\!\!\!\!\hspace{-0.5cm}\left(\sigma \displaystyle\sum_{i\in \mathcal{G}\cup \mathcal{G}^+\cup \mathcal{G^{-}} \cup\mathcal{G}_{g}^{+};\forall g}\hspace{-0.2cm}c^{\rm G}_{i}g^{(t)}_{i}+\right. \nonumber \\
&\quad\quad\quad\quad\quad\left.+\sigma \displaystyle\sum_{j\in \mathcal{D}}c^{\rm S}_{j}r^{(t)}_{j}\right); \label{OPC}
\end{align}
subject to
\begin{align}
\sum_{i \in \Psi_{n}^{\rm G}}g^{(t)}_{i}&-\sum_{k \mid o(k)=n}f^{(t)}_{k}+\sum_{k \mid r(k)=n}f^{(t)}_{k} +\sum_{j \in \Psi_{n}^{\rm D}}r^{(t)}_{j} \nonumber\\
 & = \sum_{j\in \Psi_{n}^{\rm D}}p^{(t)}_{j}: \lambda_{n}^{(t)};\; \forall n \in \mathcal{N}; \forall t \in \mathcal{T}\label{balance}\\
f^{(t)}_{k}=b_{k} \tilde{x}^{(t)}_{k} &(\theta_{o(k)}^{(t)}-\theta_{r(k)}^{(t)}): \phi^{(t)}_{k};\; \forall k \in \mathcal{L}\cup \mathcal{L}^+;\forall t \in \mathcal{T} \label{flow}\\
\theta^{(t)}_{n}&=0: \chi^{(t)}_{n};\; n:\text{slack}; \forall t \in \mathcal{T} \label{refer angle}\\
f^{(t)}_{k} &\leq f_{k}^{\rm max}: \hat{\phi}^{(t)}_{k};\; \forall k \in \mathcal{L}\cup \mathcal{L}^+; \forall t \in \mathcal{T}\label{flow upper limit}\\
f^{(t)}_{k}&\geq -f_{k}^{\rm max}: \check{\phi}^{(t)}_{k};\; \forall k \in \mathcal{L}\cup \mathcal{L}^+; \forall t \in \mathcal{T}\label{flow lower limit}\\
\theta_{n}^{(t)} & \leq \pi: \hat{\xi}^{(t)}_{n};\; \forall n \in \mathcal{N}\backslash n:\text{slack}, \forall t \in \mathcal{T} \label{angle upper limit}\\
\theta_{n}^{(t)} & \geq -\pi: \check{\xi}^{(t)}_{n};\; \forall n \in \mathcal{N}\backslash n:\text{slack}, \forall t \in \mathcal{T} \label{angle lower limit}\\
g^{(t)}_{i} &\geq 0;\; \nonumber\\
&\forall i\in \mathcal{G} \cup \mathcal{G}^{+}\cup \mathcal{G}_{g}^{+} \cup \mathcal{G}^{-};\forall g ; \forall t \in \mathcal{T} \label{gener lower limit}\\
r^{(t)}_{j}& \geq 0;\; \forall j \in \mathcal{D}, \forall t \in \mathcal{T}\label{loadshed lower limit}\\
%
p^{(t)}_{j} &= u_{j}^{{\rm D}(t)}: \alpha_{j}^{{\rm D}(t)};\; \forall j \in \mathcal{D}, \forall t \in \mathcal{T}\label{demand upper limit}
\end{align}
\begin{align}
g^{(t)}_{i} \leq &u_{i}^{{\rm G}(t)}\tilde{y}_i^{(t)}: \varphi_{i}^{{\rm G}(t)};\nonumber\\
\forall i \in &\mathcal{G} \cup \mathcal{G}^{+} \cup \mathcal{G}^{-} \cup \mathcal{G}_{g}^{+}; \forall g; \forall t \in \mathcal{T}  \label{gener upper limit}\\
r^{(t)}_{j} \leq &e^{(t)}_j u_{j}^{{\rm D}(t)}: \varphi_{j}^{{\rm D}(t)};\; \forall j\in \mathcal{D}, \forall t \in \mathcal{T} \label{loadshed upper limit}\\
%
%
u^{{\rm G}(t)}_{i} =  \bar u^{\rm G}_{i}&-\hat u^{\rm G}_{i}z^{{\rm G}(t)}_{i};\nonumber\\
\forall i\in \mathcal{G}& \cup \mathcal{G}^{+}\cup \mathcal{G}^{-}\cup \mathcal{G}_{g}^{+} ;\forall g;\forall t \in \mathcal{T} \label{uncer1}\\
u^{{\rm D}{(t)}}_{j} =  \bar u^{{\rm D}}_{j}&h_{\mu,j}^{(t)}+\hat u^{{\rm D}}_{j}h_{\sigma,j}^{(t)}z^{{\rm D}(t)}_{j};\nonumber\\
\forall j\in \mathcal{D}&;\forall t \in \mathcal{T} \label{uncer2}\\
\sum_{i\in \mathcal{G} \cup \mathcal{G}^{+}\cup \mathcal{G}^{-}\cup \mathcal{G}_{g}^{+} ;\forall g} &z^{{\rm G}(t)}_{i} \le \Gamma^{\rm G}\left( \tilde{y}^{(t)}_{i},\forall i \right);\;\forall t   \in \mathcal{T} \label{uncer3}\\
\sum_{j\in \mathcal{D}} z^{{\rm D}(t)}_{j} &\le \Gamma^{\rm D};\;\forall t \in \mathcal{T}    \label{uncer4}\\
z^{{\rm G}(t)}_{i} \in \{0,1\};\forall i\in &\mathcal{G} \cup \mathcal{G}^{+} \cup \mathcal{G}^{-}\cup \mathcal{G}_{g}^{+};\forall g;\forall t \in \mathcal{T} \label{uncer5}\\
z^{{\rm D}(t)}_{j} \in \{0,1\};\forall j\in &\mathcal{D};\forall t \in \mathcal{T} \label{uncer6}
\end{align}
\begin{align}
z^{{\rm G}(t)}_{i}  \le  \tilde{y}_i^{(t)};\forall i\in &\mathcal{G} \cup \mathcal{G}^{+}\cup \mathcal{G}^{-}\cup \mathcal{G}_{g}^{+} ;\forall g;\forall t \in \mathcal{T} \label{uncer7}
\end{align}

Equation (\ref{OPC}) represents the worst operational costs, which maximizes generation and load-shedding costs.
Constraints (\ref{balance})-(\ref{loadshed lower limit}) represent operational constraints such as setting power balance, line flows, reference bus, flow and voltage angle limits, etc. Check reference \cite{MinguezGandBertrand:16} for more details about these constraints. Restriction~(\ref{demand upper limit}) makes the level of demand match the uncertain demand variable.
Constraint~(\ref{gener upper limit}) is novel and sets the power generation to be lower than the uncertain generation capacity variable multiplied by the binary variable $\tilde{y}_i^{(t)}$, which establishes if the corresponding generator is active for period $t$. In case it is not active, i.e. $\tilde{y}_i^{(t)}=0$, the power generation is set to zero. Constraint (\ref{loadshed upper limit}) limits load-shedding to a percentage of the uncertain demand variable.

Constraints (\ref{uncer1})-(\ref{uncer6}) define the polyhedral uncertainty set analogously as it is done in \cite{MinguezGandBertrand:16}. Random generation capacity $u^{G(t)}_i$ depends on binary variable $z^{{\rm G}(t)}_{i}$, if the binary variable $z^{{\rm G}(t)}_{i}$ is 1, maximum generation capacity is set to the nominal value $\bar{u}^{G(t)}_i$ minus the maximum deviation allowed from the nominal value $\hat{u}^{G(t)}_i$, otherwise maximum generation capacity is set to the nominal value $\bar{u}^{G(t)}_i$.
 Analogously with uncertain demands $u^{{\rm D}{(t)}}_{j}$, although in this particular case, demand nominal values and dispersion are allowed to evolve during the time horizon using parameters $h_{\mu,j}^{(t)}$ and $h_{\sigma,j}^{(t)}$. These parameters allow to introduce the possible evolution of demands and their uncertainties (see reference \cite{MinguezGandBertrand:16} for more details).
 The level of uncertainty is controlled throughout the uncertainty budgets $\Gamma^{\rm G}$ and $\Gamma^{\rm D}$, which sets the maximum number of generators whose maximum capacity might be different from their nominal values and the maximum load levels that might change with respect to nominal values, respectively. However, unlike in reference \cite{MinguezGandBertrand:16} where the generation uncertainty budget was constant, in this case the uncertainty budget for each time period $\Gamma^{\rm G}$ is a function of the number of active generators for each period, i.e. $\tilde{y}^{(t)}_{i}, \forall i$. Let remind the reader that the uncertainty budget is the maximum number of generators whose maximum capacity is allowed to depart from their nominal values, if the number of generator increases, the uncertainty budget should increase to keep an analogous level of protection against uncertainty. Note that what is the appropriate selection of this function $\Gamma^{\rm G}\left(\tilde{y}^{(t)}_{i}, \forall i \right)$ is out of the scope of the paper.
 Finally, constraint  (\ref{uncer7}) is also novel for this work and sets the binary variables related to generators to zero in case generators are not active at time period $t$, thus they can not account for uncertainty budget in (\ref{uncer3}).

\vspace{-0.3cm}
\section{Proposed Decomposition Method}\label{Descomposicion}
\vspace{-0.1cm}
%
The aim of this section is to extend the solution procedure presented in \cite{MinguezGandBertrand:16} to solve the robust DTNRGEP problem described in Section~\ref{Model}. Since the decomposition method has a bi-level structure, the first step is to merge the initial three-level formulation (\ref{eq1})-(\ref{uncer7}) into a two-level problem.
\subsection{Second-level formulation: subproblem}
For given values for the first-stage variables $\tilde{x}_k^{(t)}$ and $\tilde{y}_{i}^{(t)}$ for each time period, the problem set out by (\ref{OPC})-(\ref{uncer7}) might be decomposed into the following single-level maximization problem for each time period $t;\;\forall t\in \mathcal{T}$:
\begin{align}
&c_{\rm op}^{(t)}=\!\!\!\!\mm{Maximize}{{\bfg u},\lambda_{n}^{(t)},\phi^{(t)}_{k},\chi^{(t)}_{n},\hat{\phi}^{(t)}_{k},\\ \check{\phi}^{(t)}_{k},\hat{\xi}^{(t)}_{n},\check{\xi}^{(t)}_{n},\alpha_{j}^{{\rm D}(t)},\varphi_{i}^{{\rm G}(t)},\varphi_{j}^{{\rm D}(t)}}{}\nonumber\\
\vspace{0.2cm}
&{\left\{\!
\begin{array}{c}
 \displaystyle \sum_{k \in \mathcal{L}}\Bigl(\hat{\phi}^{(t)}_{k}-\check{\phi}^{(t)}_{k}\Bigr)f_{k}^{\rm max}+\!\!\!\displaystyle\sum_{n\in \mathcal{N}\backslash n:\textrm{slack}}\pi\Bigl(\hat{\xi}^{(t)}_{n}-\check{\xi}^{(t)}_{n}\Bigr) \\
 +\!\!\!\displaystyle\sum_{i\in \mathcal{G} \cup \mathcal{G}^{+} \cup \mathcal{G}^{-}\cup \mathcal{G}_{g}^{+};\forall g}\!\!\Bigl(u_{i}^{{\rm G}(t)}\tilde{y}_i^{(t)}\varphi_{i}^{{\rm G}(t)}\Bigr)\\
 +\!\!\!\displaystyle\sum_{j\in \mathcal{D}}\!\!\Bigl(u_{j}^{{\rm D}(t)}\alpha_{j}^{{\rm D}(t)}\!\!+\!\!e^{(t)}_j u_{j}^{{\rm D}(t)}\varphi_{j}^{{\rm D}(t)}\!\! \Bigr)
\end{array}\right\}
}\nonumber\\
\label{subproblem1}
\end{align}

subject to:
\begin{align}
&\lambda_{n(i)}^{(t)}+\varphi_{i}^{{\rm G}(t)} \leq  \fraca{\sigma}{(1+I)}c^{\rm G}_{i}; \nonumber\\
&\forall i\in \mathcal{G} \cup \mathcal{G}^{+} \cup \mathcal{G}^{-}\cup \mathcal{G}_{g}^{+};\forall g; \forall t \in \mathcal{T} \label{dual blockgen}\\
&-\lambda_{n(j)}^{(t)}+\alpha_{j}^{{\rm D}(t)}  \leq  0; \; \forall j \in \mathcal{D};\forall t \in \mathcal{T} \label{dual blockdemand}\\
&\lambda_{n(j)}^{(t)}+\varphi_{j}^{{\rm D}(t)}  \leq   \fraca{\sigma}{(1+I)}c^{\rm S}_{j}; \;\forall j\in \mathcal{D};\forall t \in \mathcal{T} \label{dual loadshed}\\
&-\lambda_{o(k)}^{(t)}+\lambda_{r(k)}^{(t)}+\phi^{(t)}_{k}+\hat{\phi}^{(t)}_{k}+\check{\phi}^{(t)}_{k} =  0; \nonumber\\
& \forall k\in \mathcal{L}\cup\mathcal{L}^+;\forall t \in \mathcal{T} \label{dual flow}\\
%
&-\sum_{k\mid o(k)=n}b_{k}\tilde{x}^{(t)}_{k}\phi^{(t)}_{k}+\sum_{k\mid r(k)=n}b_{k}\tilde{x}^{(t)}_{k}\phi^{(t)}_{k}    \nonumber \\
&+\hat{\xi}^{(t)}_{n}+\check{\xi}^{(t)}_{n}  =  0 ; \;\forall n \in \mathcal{N}\backslash n:\text{slack};\forall t \in \mathcal{T} \label{dual angle}\\
&-\sum_{k\mid o(k)=n}b_{k}\tilde{x}^{(t)}_{k}\phi^{(t)}_{k}+\sum_{k\mid r(k)=n}b_{k}\tilde{x}^{(t)}_{k}\phi^{(t)}_{k}    \nonumber \\
&+\chi^{(t)}_{n} = 0;\;n:\text{slack};\forall t \in \mathcal{T} \label{dual anglerefer}\\
&-\infty \leq  \lambda_{n}^{(t)}  \leq\infty; \;\forall n \in \mathcal{N};\forall t \in \mathcal{T} \label{dual lambda}\\
&-\infty \leq \phi^{(t)}_{k}  \leq\infty; \;\forall k \in \mathcal{L}\cup\mathcal{L}^+;\forall t \in \mathcal{T} \label{dual phik}\\
&-\infty \leq \chi^{(t)}_{n} \leq\infty; \; n:\text{slack};\forall t \in \mathcal{T} \label{dual chi}\\
& \hat{\phi}^{(t)}_{k}\leq 0; \;\forall k  \in \mathcal{L}\cup\mathcal{L}^+;\forall t \in \mathcal{T}\label{dual phikmax}\\
& \check{\phi}^{(t)}_{k} \geq  0; \;\forall k \in \mathcal{L}\cup\mathcal{L}^+;\forall t \in \mathcal{T}\label{dual phikmin}\\
& \hat{\xi}^{(t)}_{n}\leq 0; \;\forall n  \in \mathcal{N}\backslash n:\text{slack};\forall t \in \mathcal{T}\label{dual ximax}\\
& \check{\xi}^{(t)}_{n} \geq  0; \; \forall n \in \mathcal{N}\backslash n:\text{slack};\forall t \in \mathcal{T}\label{dual ximin}\\
&-\infty \leq  \alpha_{j}^{{\rm D}(t)}  \leq\infty; \; \forall j \in \mathcal{D};\forall t \in \mathcal{T} \label{dual alpha}\\
& \varphi_{i}^{{\rm G}(t)} \leq 0; \; \forall i  \in \mathcal{G} \cup \mathcal{G}^{+} \cup \mathcal{G}^{-}\cup \mathcal{G}_{g}^{+};\forall g;\forall t \in \mathcal{T}\label{dual varphiG}\\
& \varphi_{j}^{{\rm D}(t)} \leq 0; \; \forall j  \in \mathcal{D};\forall t \in \mathcal{T}\label{dual varphiD}\\
& \mbox{Constraints (\ref{uncer1})-(\ref{uncer7})}\label{uncerb}
\end{align}

Subproblems (\ref{subproblem1})-(\ref{uncerb}) result from substituting in problem (\ref{OPC})-(\ref{uncer7}) for each time period $t$, the third-level problem by its dual. An important aspect for the resolution of subproblems is the linearization of bilinear terms in (\ref{subproblem1}), i.e., $\sum_{i\in \mathcal{G} \cup \mathcal{G}^{+} \cup \mathcal{G}^{-}\cup \mathcal{G}_{g}^{+};\forall g}\Bigl(u_{i}^{{\rm G}(t)}\tilde{y}_i^{(t)}\varphi_{i}^{{\rm G}(t)}\Bigr)+\sum_{j\in \mathcal{D}}\Bigl(u_{j}^{{\rm D}(t)}\alpha_{j}^{{\rm D}(t)}+e^{(t)}_j u_{j}^{{\rm D}(t)}\varphi_{j}^{{\rm D}(t)}\Bigr)$. The linearization process is described in more detail in \cite{MinguezG:15}. Note that variable $\tilde{y}_i^{(t)}$ is considered a parameter within our subproblem.
These subproblems provide the uncertain parameter values ${\bfg u}^{(t)}$ within the uncertainty sets to give the least desirable operational costs for each year.

The resulting formulation associated with subproblems is a mixed-integer linear programming problem, which can be solved efficiently by using state-of-the-art mixed-integer mathematical programming solvers such as CPLEX or Gurobi.
\vspace{-0.2cm}
\subsection{First-level formulation: master problem}
First level formulation corresponds to the master problem. Thus, for given values of the uncertain parameter values ${u}^{(t)}\in \mathcal{U}^{(t)}$ for each year obtained from second-level subproblems, the master problem at iteration $\nu$ corresponds to:
\begin{align}
  \mm{Minimize}{\begin{array}{c}
                  x_k^{(t)},\tilde{x}_k^{(t)}, y_k^{(t)}, \tilde{y}_k^{(t)}, \\
                  g^{(t)}_{i,l}, p^{(t)}_{j,l}, r^{(t)}_{j,l}, \theta^{(t)}_{n,l}, f^{(t)}_{k,l} \\
                  l= 1,\ldots,\nu-1
                \end{array}}{\displaystyle\sum_{t \in \mathcal{T}} \fraca{1}{(1+I)^{t-1}}\left(\displaystyle\sum_{k \in \mathcal{L}^+}c_k^{\rm LI} x_k^{(t)}\right. \nonumber\\
\left.+\displaystyle\sum_{i \in \mathcal{G}^{+}\cup \mathcal{G}_{g}^{+};\forall g}c_i^{\rm GI} y_i^{(t)}+\gamma^{(t)}\right);} \nonumber\\
\label{master5}
\end{align}
subject to
\begin{align}
\gamma^{(t)}   \ge \fraca{\sigma}{(1+I)} \hspace{-0.2cm}&\displaystyle\sum_{\begin{array}{c} i\in \mathcal{G} \cup \mathcal{G}^{+}\\ \cup \mathcal{G}^{-} \cup \mathcal{G}_{g}^{+} ;\forall g\end{array}} \hspace{-0.2cm} c^{\rm G}_{i}g^{(t)}_{i,l}+ \fraca{\sigma}{(1+I)} \displaystyle\sum_{j\in \mathcal{D}}c^{\rm S}_{j}r^{(t)}_{j,l};\;\nonumber\\
&\quad \quad\quad\quad\forall t \in \mathcal{T}, l = 1,\ldots,\nu-1\label{master7}\\
\gamma^{(t)}& \ge 0;\;\forall t\in \mathcal{T}\\
\textrm{Constraints } &(\ref{eq2})-(\ref{eq2n})\label{masterinicial}\\
\textrm{Equations } &(\ref{balance})-(\ref{loadshed upper limit});\;l = 1,\ldots,\nu-1. \label{masterfinal}
\end{align}

Constraint (\ref{master7}) are primal decomposition cuts, while variables $\gamma^{(t)}$ relates to year on year operational costs. The master problem  includes one variable $g^{(t)}_{i,l}$, $p^{(t)}_{j,l}$, $r^{(t)}_{j,l}$, $\theta^{(t)}_{n,l}$ and $f^{(t)}_{k,l}$ for each year and for each realization of the uncertain parameters obtained from the subproblem (\ref{subproblem1})-(\ref{uncerb}) at every iteration.

\vspace{-0.2cm}
\subsection{Solution method}
The solution method consists of iteratively solving the following master and subproblems at each iteration $\nu$:
\begin{itemize}
    \item {\bf Master problem:} For given realizations of the uncertain parameters obtained from subproblems at the previous iterations, new first-stage variables $\tilde{x}_{k,\nu}^{(t)}$, $\tilde{y}_{i,\nu}^{(t)}$ values are calculated by means of (\ref{master5})-(\ref{masterfinal}). The optimal objective function lower bound is updated $z^{(\rm lo)}=\sum_{t \in \mathcal{T}} \frac{1}{(1+I)^{t-1}}(\sum_{k \in \mathcal{L}^+}c_k^{\rm LI} x_{k,\nu}^{(t)}+\sum_{i \in \mathcal{G}^{+}\cup \mathcal{G}_{g}^{+};\forall g}c_i^{\rm GI} y_{i,\nu}^{(t)}+\gamma^{(t)})$.

  \item {\bf Subproblems, one for each year:} For given values of the first-stage decision variables $\tilde{x}_{k,\nu}^{(t)}$, $\tilde{y}_{i,\nu}^{(t)}$, uncertain parameters within the uncertainty set giving the least desirable operational costs (\ref{OPC}), i.e. ${\bfg u}^{(t)}_{\nu}$ and $c_{{\rm op},\nu}^{(t)}$, respectively, are calculated by solving  subproblems in (\ref{subproblem1})-(\ref{uncerb}). The optimal objective function upper bound is updated $z^{(\rm up)}=\sum_{t\in \mathcal{T}} \frac{1}{(1+I)^{t-1}}(\sum_{k \in \mathcal{L}^+}c_k^{\rm LI} x_{k,\nu}^{(t)}+\sum_{i \in \mathcal{G}^{+}\cup \mathcal{G}_{g}^{+};\forall g}c_i^{\rm GI} y_{i,\nu}^{(t)}+c_{{\rm op},\nu}^{(t)})$.
\end{itemize}

This iterative process is repeated until the absolute value  of the relative difference between upper and lower bounds is below a selected threshold. For more details about the structure of the algorithm, applied just for the capacity expansion problem, check reference \cite{MinguezGandBertrand:16}.

The advantage of this bi-level formulation given by (\ref{subproblem1})-(\ref{uncerb})  and (\ref{master5})-(\ref{masterfinal}) is that it has the same problem structure than that defined by \cite{ZengZ:13}, and therefore the proposed column-and-constraint generation method guarantees convergence to a global optimum.

\vspace{-0.1cm}
\section{Examples}
\label{s4}

In this section, the numerical results for an illustrative example based on Garver system \cite{Garver:70} and a case study using the IEEE $118$-bus test system \cite{118system} are portrayed to analyze the joint study of transmission network and generation expansion planning.

All numerical tests have been implemented and solved using CPLEX within GAMS \cite{brooke1998gams} on a Windows DELL PowerEdge R$920$ server with two Intel Xeon E$7-4820$ processors clocking at $2$ GHz and $768$ Gb of RAM. The stopping tolerance for all cases is equal to $\varepsilon=10^{-6}$.
\vspace{-0.2cm}
\subsection{Illustrative Example. Garver's 6-bus System}
\begin{figure}[H]
  \begin{center}
  \includegraphics[width=0.44\textwidth]{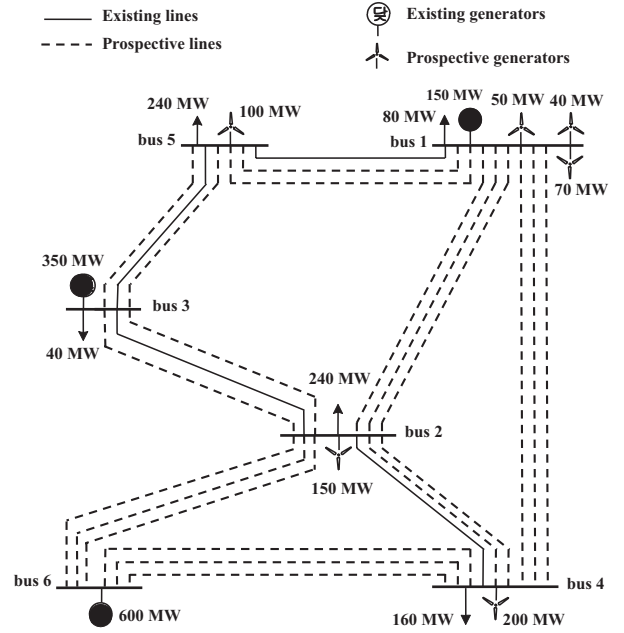}
    \caption{Garver's 6-bus test system.}
    \label{Figure1}
  \end{center}
\end{figure}

The model is initially tested in the Garver 6-bus system depicted in Figure~\ref{Figure1}. This system is composed of 6 buses, 3 generators, 5 levels of inelastic demand and 6 lines. Data for generation and demand capacities, and supply and bidding prices are given in \cite{MinguezG:15}. The load-shedding cost is equal to hundred times the bidding price for each level of demand. Line data are obtained from Table I of reference \cite{GarcesCGR:09}, including construction costs.

Regarding expansion possibilities, 6 generation units and 3 transmission lines between each pair of buses could be installed. The characteristics of possible generation units, which are assumed to be wind units, are shown in Table~\ref{NewGen}. The first three generation units belong to one group ($g=1$) to be constructed and installed in phases (at different times) at the same bus 1, so in case the model decides to include them in the generation expansion planning they must be installed sequentially, i.e. $\mathcal{G}_{1}^{+}\equiv\{4,5,6\}$. In addition, it is known that the existing generator at node 1 is going to be dismantled at time period 8 because it reaches its useful life, i.e. $\mathcal{G}^{-}\equiv\{1\}$. The rest of sets associated with generation are $\mathcal{G}\equiv\{2,3\}$ and $\mathcal{G}^{+}\equiv\{7,8,9\}$.
Operational and investment costs for generators presented in Table~\ref{NewGen} come from \cite{de2013optimal}.
In what regards to the characteristics of the new possible transmission lines, they are also attained from Table I of reference~\cite{GarcesCGR:09}. The maximum available investment budget for transmission lines is 40 million euros.

The time horizon considered is 25 years and the discount rate is 10\%. The weighted factor $\sigma$ is equal to the number of hours in one year, i.e. 8760, so that the load-shedding and power generation costs are related to years, which can be comparable with the annualized investment cost.

\vspace{-0.2 cm}

\begin{table}[H]
    \caption{Candidate renewable generators related to Garver's 6-bus test system.}
    \label{NewGen}
    \centering
     \renewcommand{\tabcolsep}{1mm}
\begin{tabular}{ccccc}
\hline
\multirow{2}{*}{Gen.} & \multirow{2}{*}{Bus} &  Power Capacity & O\&M Costs       & Investment cost  \\
     &                &   $\bar u^{\rm G}_{i}$ (MW)      & (\euro/MWh)& (M\euro) \\
\hline
4 & 1 & $50$ & $17.8$ & $50$ \\\hline
5 & 1 & $70$ & $17.5$ & $80$ \\\hline
6 & 1 & $40$ & $17.5$ & $40$\\\hline
7 & 2 & $150$ & $16.5$ & $200$\\\hline
8 & 4 & $200$ & $15.0$ & $198$\\\hline
9 & 5 & $100$ & $17.0$ & $110$\\\hline
\end{tabular}
\end{table}

\vspace{-0.2 cm}

Regarding the uncertainty sets, the maximum capacity of conventional generators can decrease a maximum of 50\% with respect to their nominal values, i.e.  $\hat{u}^{\rm G}_{i}=0.5\bar u^{\rm G}_{i};\;i=1,2,3$, while for renewable generators their maximum capacity can decrease a 100\% with respect to their nominal values, i.e.  $\hat{u}^{\rm G}_{i}=\bar u^{\rm G}_{i};\;i=4,\ldots,9$. Load levels may change a maximum of 20\% with respect to their nominal values, i.e. $\hat{u}^{\rm D}_{j}=0.2\bar u^{\rm D}_{j};\;j=1,\ldots,5$. Finally, annual growth rates for load nominal values and dispersion are equal to 1.2\%, i.e. $h_{\mu,j}^{(t)}=h_{\sigma,j}^{(t)}=1.012^{(t-1)}$.

Four case studies are analyzed, considering two different combinations of uncertainty budgets associated with generation capacities and demands and two different generation investment budgets. It is worth stressing that the inclusion of new generators implies updating the uncertainty budget associated with generation, we use the following step function to define $\Gamma^{\rm G}\left(\tilde{y}^{(t)}_{i}, \forall i \right)$: if 1 or 2 generators are built, the generation uncertainty budget increases in one unit; if 3 or 4 generators are built, the uncertainty budget increases in two units; and finally, if 4 or 5 generators are built, the uncertainty budget increases in three units.

Results about investment cost, and lines and generators built for each case study using the proposed model are given in Table~\ref{GarverSolution}.
From this table the following observations are pertinent:

\begin{itemize}
\item[{\bf{$\bullet$}}] For the first combination of uncertainty budgets, $\Gamma^{G}$ and $\Gamma^{D}$ equal to 1 and 2, lines and generators built are the same. Note that the candidate generator that remains to be built has an investment cost of 200 million euros, thus the rise considered in the generation investment budget does not allow the construction of the remaining generator and no more lines are needed.

\item[{\bf{$\bullet$}}] Total cost for cases a) and c) is 14627.918 and 64927.727 million euros, respectively. Therefore, the higher uncertainty is, the higher total cost is. Note that lines are built in earlier periods when the total cost is higher.

\item[{\bf{$\bullet$}}] Total cost for case studies b) and d) is 14627.918 and 63081.771 million euros, respectively. Note that the configuration of built generators changes and lines are built early in d) due to the higher cost.

\item[{\bf{$\bullet$}}] In case studies c) and d), the increase in investment budget for generators allows to build one more line or more generators.

\item[{\bf{$\bullet$}}] It is worth stressing the complexity of this problem. For the same uncertainty budgets, the solutions associated with different investment budgets are not incremental. Thus the necessity to consider the year-by-year dynamic of the problem.

\end{itemize}

Regarding computational tractability, the number of iteration required are 16, 15, 10 and 16, respectively, for study cases a), b), c) and d). The maximum computing time is eight hours for the case a).

\subsection{IEEE 118-bus test system example}
We also apply the proposed model on a bigger and more realistic case using the IEEE 118-bus test system~\cite{118system}.
The system is composed of 118 buses, 186 existing lines, 54 generators and 91 loads. We assume that the generator located at bus 4 stops working at period 8 because it reaches its useful life, , i.e. $\mathcal{G}^{-}\equiv\{4\}$. Generation capacities and demand loads can be found in \cite{MinguezG:15}. The load-shedding cost is ten times the bidding price of each level of demand. In addition, the same 61 existing lines given in \cite{MinguezG:15} can be duplicated to build additional lines. Data for all lines are taken from \cite{118system}. The characteristics of new possible generators are shown in Table~\ref{NewGen118} and there are two groups of generators to be installed sequentially at buses 4 and 20, i.e. $\mathcal{G}^{+}_1\equiv\{56,57,58\}$ and $\mathcal{G}^{+}_2\equiv\{64,65,66\}$.
The investment budgets for the generators and transmission lines are 1500 and 100 million euros, respectively. The discount rate is 10\% and the time horizon is 25 years.

Due to uncertainty, conventional generators can decrease a 50\% with respect to their nominal values, i.e.  $\hat{u}^{\rm G}_{i}=0.5\bar u^{\rm G}_{i};\;i=1,\ldots,54$, while for renewable generators their maximum capacity can decrease a 100\% with respect to their nominal values, i.e.  $\hat{u}^{\rm G}_{i}=\bar u^{\rm G}_{i};\;i=55,\ldots,84$. Demand levels may change a maximum of 50\%  with respect to their nominal values. Annual growth rates for load nominal values and dispersion are equal to 1.2\%, i.e. $h_{\mu,j}^{(t)}=h_{\sigma,j}^{(t)}=1.012^{(t-1)}$.
\vspace{-0.2 cm}
\begin{table}[H]
    \caption{Results for Garver's 6-bus test system illustrative example.}
    \label{GarverSolution}
    \centering
     \renewcommand{\tabcolsep}{1mm}
\begin{tabular}{cccccccc}
\hline
Case  & Input & Inv. Cost  & \multicolumn{3}{c}{New lines}  & \multicolumn{2}{c}{New generators}  \\
\hhline{~~~||---||--}
   study   & $\Gamma^{G}$, $\Gamma^{D}$, $\Pi_{G}$(M\euro) &  (M\euro)     & From & To & Period & Bus & Period \\

\hline
\multirow{6}{*}{a)}&\multirow{6}{*}{1, 2, 350}  & \multirow{6}{*}{384.802} & 1 & 5 & 5 & 4 & 8  \\
                                                                               & & & 2 & 3 & 1 & 5 & 3  \\
                                                                               & & & 2 & 6 & 1 & 1 & 1  \\
                                                                               & & & 2 & 6 & 8 & 1 & 2   \\
                                                                               & & & 2 & 6 & 1 & 1 & 3   \\
                                                                               & & & 3 & 5 & 1 & - & -  \\
                                                                               & & & 4 & 6 & 1 & - & -  \\
                                                                               & & & 4 & 6 & 1 & - & -  \\
\hline
\multirow{6}{*}{b)}&\multirow{6}{*}{1, 2, 450} & \multirow{6}{*}{384.802} & 1 & 5 & 5 & 4 & 8  \\
                                                                               & & & 2 & 3 & 1 & 5 & 3  \\
                                                                               & & & 2 & 6 & 1 & 1 & 1  \\
                                                                               & & & 2 & 6 & 8 & 1 & 2   \\
                                                                               & & & 2 & 6 & 1 & 1 & 3   \\
                                                                               & & & 3 & 5 & 1 & - & -  \\
                                                                               & & & 4 & 6 & 1 & - & -  \\
                                                                               & & & 4 & 6 & 1 & - & -  \\

\hline
\multirow{6}{*}{c)}&\multirow{6}{*}{2, 3, 350}  & \multirow{6}{*}{385.190} & 1 & 5 & 1 & 4 & 8  \\
                                                                               & & & 2 & 6 & 2 & 5 & 1  \\
                                                                               & & & 2 & 6 & 1 & 1 & 1  \\
                                                                               & & & 2 & 6 & 1 & 1 & 4  \\
                                                                               & & & 3 & 5 & 1 & 1 & 5  \\
                                                                               & & & 4 & 6 & 1 & - & -  \\
                                                                               & & & 4 & 6 & 1 & - & -  \\
\hline
\multirow{6}{*}{d)}&\multirow{6}{*}{2, 3, 450}  & \multirow{6}{*}{475.124} & 1 & 5 & 1 & 2 & 2  \\
                                                                               & & & 2 & 3 & 2 & 5 & 1  \\
                                                                               & & & 2 & 6 & 6 & 1 & 1  \\
                                                                               & & & 2 & 6 & 1 & 1 & 2 \\
                                                                               & & & 2 & 6 & 1 & 1 & 7  \\
                                                                               & & & 3 & 5 & 1 & - & -  \\                                                                               & & & 4 & 6 & 1 & - & -  \\
                                                                               & & & 4 & 6 & 1 & - & -  \\
\hline
\end{tabular}
\end{table}

\vspace{-0.3 cm}

\begin{table}[H]
    \caption{Candidate generators related to IEEE 118-bus test system example.}
    \label{NewGen118}
    \centering
     \renewcommand{\tabcolsep}{1mm}
\begin{tabular}{ccccc}
\hline
\multirow{2}{*}{Number} & Bus &  Power Capacity & O\&M Cost       & Investment cost  \\
                   &  &   (MW)      & (\euro/MWh)& (M\euro) \\
\hline
55 & 1 & $90$ & $15.2$ & $120$ \\\hline
56 & 4 & $50$ & $17.8$ & $50$ \\\hline
57 & 4 & $70$ & $17.5$ & $80$\\\hline
58 & 4 & $40$ & $17.5$ & $40$\\\hline
59 & 6 & $100$ & $17.1$ & $110$\\\hline
60 & 10 & $180$ & $15.3$ & $184$ \\\hline
61 & 14 & $100$ & $17.0$ & $110$\\\hline
62 & 14 & $90$ & $15.2$ & $120$\\\hline
63 & 18 & $150$ & $16.0$ & $145$ \\\hline
64 & 20 & $50$ & $17.6$ & $50$ \\\hline
65 & 20 & $50$ & $17.6$ & $50$\\\hline
66 & 20 & $60$ & $15.4$ & $55$\\\hline
67 & 21 & $130$ & $16.5$ & $135$\\\hline
68 & 22 & $200$ & $15.0$ & $198$\\\hline
69 & 27 & $80$ & $15.9$ & $90$\\\hline
70 & 38 & $110$ & $16.7$ & $123$ \\\hline
71 & 39 & $200$ & $15.1$ & $200$ \\\hline
72 & 50 & $90$ & $17.0$ & $118$\\\hline
73 & 51 & $150$ & $16.6$ & $153$\\\hline
74 & 62 & $110$ & $16.8$ & $103$\\\hline
75 & 75 & $110$ & $16.9$ & $147$ \\\hline
76 & 80 & $170$ & $16.0$ & $164$\\\hline
77 & 88 & $200$ & $14.9$ & $198$\\\hline
78 & 93 & $100$ & $17.0$ & $110$ \\\hline
79 & 94 & $200$ & $15.0$ & $198$ \\\hline
80 & 96 & $140$ & $16.3$ & $158$ \\\hline
81 & 101 & $170$ & $15.2$ & $180$\\\hline
82 & 114 & $190$ & $15.5$ & $191$\\\hline
83 & 116 & $110$ & $16.6$ & $112$\\\hline
84 & 118 & $90$ & $17.3$ & $102$ \\\hline
\end{tabular}
\end{table}

Using the following uncertainty budgets $\Gamma^{G}=15$ and $\Gamma^{D}=20$, the DTNRGEP approach provides an investment cost of 1599.574 million euros for constructing the generators shown in Table~\ref{LinesGen} and the lines 187, 189, 191, 192, 203, 204, 205, 206, 207, 211, 223, 226, 241. Lines 189 and 204 are constructed at period 15 and 3, respectively, while the rest of lines are constructed at the beginning of the time period considered. In terms of operational costs, a total of 1.069 billion euros are needed, of which 360.961 million euros belong to load-shedding.

\begin{table}[H]
    \caption{New generators for IEEE 118-bus test system example.}
    \label{LinesGen}
    \centering
     \renewcommand{\tabcolsep}{1mm}
\begin{tabular}{cccccc}
\hline
Number & Bus & Period & Number & Bus & Period \\
\hline
56 & $4$ & $22$  & 63 & $18$ & $1$  \\
64 & $20$ & $1$ & 65 & $20$ & $2$ \\
66 & $20$ & $3$ & 71 & $39$ & $1$ \\
73 & $51$ & $1$  & 74 & $62$ & $1$ \\
76 & $80$ & $1$ & 77 & $88$ & $1$ \\
79 & $94$ & $1$  & 82 & $114$ & $1$ \\
\hline
\end{tabular}
\end{table}

\vspace{-0.1 cm}

Note that as in the previous example, the sequential installation of the corresponding generators is respected. Thus, three generators are built at periods 1, 2 and 3, respectively, at bus 20.

Regarding computational tractability, 5 iterations are required to reach convergence in a computational time of 4 hours and 28 minutes.
\vspace{-0.2cm}
\section{Conclusions}\label{s5}
In this paper the use of robust optimization for solving the dynamic transmission and renewable generation expansion planning problem has been extended. The model put forward herein provides the initial design and the expansion plan as regards forthcoming years in terms of where and when new lines and/or generators have to be constructed. It is possible to assume that the probability distributions for the random variables (uncertainty sets) change between consecutive years. The proposed model provides an integrated approach reaching the global optimal solution, and overcomes the size limitations and computational intractability associated with this type of problem for realistic cases.

The proposed model constitutes a valuable tool for the difficult task of transforming conventional power systems into sustainable systems.
\vspace{-0.2cm}

\ifCLASSOPTIONcaptionsoff
  \newpage
\fi




\end{document}